\documentclass{article}
\usepackage{amssymb}
\begin{document}
\newcommand{\uX}{\underline{X}}
\newcommand{\uy}{\underline{y}}
\newcommand{\Aut}{{\rm Aut}}
\newcommand{\cU}{{\cal U}}
\newcommand{\cY}{{\cal X}}
\newcommand{\Fut}{{\rm Futaki}}
\newcommand{\Ric}{{\rm Ric}}
\newcommand{\oK}{$\overline{K}$}
\newcommand{\cX}{{\cal X}}
\newcommand{\cW}{{\cal W}}
\newcommand{\op}{{\rm op}}
\newcommand{\bC}{{\bf C}}
\newcommand{\bP}{{\bf P}}
\newcommand{\cS}{{\cal S}}
\newtheorem{thm}{Theorem}
\newtheorem{lem}{Lemma}
\newtheorem{cor}{Corollary}
\newtheorem{prop}{Proposition}
\newcommand{\omu}{\overline{\mu}}
\newcommand{\Ch}{{\rm Ch}}
\newcommand{\Vol}{{\rm Vol}}
\newcommand{\Tr}{{\rm Tr}}
\newcommand{\bQ}{{\bf Q}}
\newcommand{\cL}{{\cal L}}
\newcommand{\bR}{{\bf R}}
\newcommand{\cV}{{\cal V}}
\newcommand{\cO}{{\cal O}}
\newcommand{\onu}{\overline{\nu}}
\newcommand{\osigma}{\underline{\sigma}}
\newcommand{\usigma}{\underline{\sigma}}
\newcommand{\cQ}{{\cal Q}}
\newcommand{\lambdamax}{\lambda_{{\rm max}}}
\newtheorem{defn}{Definition}
\newcommand{\Hilb}{{\rm Hilb}}

\title{Stability, birational transformations and the Kahler-Einstein problem}
\author{S. K. Donaldson}
\maketitle

\section{Introduction}

This is the first in a series of papers in which we will discuss the existence problem for  Kahler-Einstein metrics on complex projective manifolds.\footnote{ An informal document {\it Discussion of the Kahler-Einstein problem}, giving some more details, is available on the webpage http:/www2.imperial.ac.uk/~skdona/.} It is well-known that, following the work of Yau  in the 1970's, the question is reduced to the \lq\lq positive'' case, so we want to decide when a Fano manifold admits a Kahler-Einstein metric. The present paper is confined to algebro-geometric aspects of the problem. We have three main purposes
\begin{itemize}
\item To state   some new definitions, of \lq\lq b-stability'' and \lq\lq \oK-stability''.
\item To indicate the application of these to the differential geometric existence problem.
\item To highlight some purely algebro-geometric questions which seem to be important in the existence problem. 
\end{itemize}

In the late 1980's, Yau suggested that the existence of Kahler-Enstein metrics should be related to  the algebro-geometric notion of {\it stability}, whose origins lie in Geometric Invariant Theory\cite{kn:Y}. This  was motivated in part by  analogy  with the Kobayshi-Hitchin correspondence  for Hermitian Yang-Mills connections on  holomorphic vector bundles. The conjecture was refined considerably by work of Tian in the 1990's {\cite{kn:T1}. We will indicate this general idea---the {\it Tian-Yau conjecture}---by the informal slogan 
$$  {\rm KE \ metric}\ \  \Longleftrightarrow \ \ \ {\rm \lq\lq algebro-geometric \ stability''}. $$
We should emphasise however that, while this slogan summarises the general idea, the exact notion of \lq\lq stability'' needs to be specified, and this specification can be considered as part of the problem.  There are a number of established notions in the literature, for example {\it asymptotic Chow stability} and {\it K-stability}: a good account of the relation between these is given  in \cite{kn:RT}. But our purpose   here is to introduce two new notions, {\it $\overline{K}$-stability} and {\it b-stability}, so that we can formulate a precise version of this conjecture which we hope will be more accessible to proof. 

Before going further some general remarks may be helpful. First, there has been a great deal of work in this area, which we will not attempt to summarise in detail here. Notably, Tian showed in \cite{kn:T1}that his notion of K-stability was a  necessary condition for the existence of a Kahler-Einstein metric, and gave an explicit example where this condition goes beyond those known before (see also (4.2) below). Second, the plethora of algebro-geometric notions of stability are all variants of the same basic idea. One expects that, among them, various definitions which are {\it a priori} different may {\it a posteriori}
turn out to be equivalent.  Since we want, ultimately, to establish a necessary and sufficient condition we can to some extent  push difficulties from one side to the other: working with a stronger (more restrictive) notion of \lq\lq algebro-geometric stability'' will, as a matter of logic, simplify the proof of ${\rm stability} \Rightarrow {\rm KE\ metric}$ while making the proof of the converse more difficult. We will explain below how this obvious remark applies in the case of \oK-stability and b-stability. The third remark is that, in the overall  existence problem for  Kahler-Einstein metrics, there is some distinction between theory and practice. One may have sufficient criteria which, while not in general necessary, can be applied to give existence statements in certain explicit cases. For example this is true of Tian's theory of the  $\alpha$-invariant. On the other hand one can envisage general, theoretical, solutions which may be difficult or impossible to apply to any specific case. This is probably the situation, as things stand at the moment, with the various notions of stability mentioned above, since they are all exceedingly difficult to verify in  examples. We will discuss this further in 4.3 below, but we note here that  our notion of \lq\lq b-stability''
comes with an integer parameter $m$---\lq\lq b-stability at multiplicity $m$''---which is related to this issue.

The general shape of all these definitions of stability is that a manifold is stable unless there is a \lq\lq destabilising object'' of an appropriate kind. The notion of b-stability is derived by extending the class of destabilising objects. This extension can be seen as part of a more general trend. In Tian's original definition of K-stability the destabilising objects were projective varieties, smooth or mildly singular, with holomorphic vector fields. In the generalisation of \cite{kn:D1} the destabilising objects were allowed to be general schemes with $\bC^{*}$ actions. In the new notion of b-stability the destabilising objects are {\it sequences} of schemes (\lq\lq webs of descendants''), related by birational transformations 
(reflected in the prefix \lq\lq b''). Another point  is that we find that we need to work not just with the  \lq\lq test configurations'' appearing in standard definitions of stability, but with more general degenerations. Thus in Section 2 we develop background and foundations to do this, review the definition of $K$-stability and state the definition  of \oK-stability. Section 4 is intended to explain the  motivation for  the definition of b-stability. We state, without complete proof, a simple existence theorem and give some general discussion. In Section 5 we collect the proofs of some auxiliary results, not used in an essential way in the core of this paper, and further examples.

The author hopes that the ideas discussed in this paper raise various questions of interest to   algebraic geometers. The main question is whether the notion of b-stability is actually needed. The definition is designed to get around a difficulty which it seems (to the author) might in principle occur in taking algebro-geometric limits of Fano manifolds under projective embeddings, defined by higher and higher powers of the anticanonical bundle. But the author does not know of any real example where this occurs (although certainly if one varies the hypotheses slightly it can: see (5.5)). It might be that, with a deeper analysis, one could show that this phenomenon does not occur (or perhaps one could show this under extra hypotheses, such as complex dimension 3). In that case one could forget the definition of b-stability and the existence theory for Kahler-Einstein metrics would be much simpler. On the other hand, it might be that one could produce actual examples of this phenomenon, and such examples  should be related to the question of whether Gromov-Haussdorf limits of Kahler-Einstein manifolds have algebraic structures.  In another direction, there is a circle of questions related to the distinction between  test configurations and general degenerations, which fall into the area of algebraic group actions. We give some simple results and examples in (5.1) and (5.2) but the author would like to understand the general picture better.

The author is very grateful to Alessio Corti, Paolo Cascini, Frances Kirwan,   Miles Reid, Julius Ross, Jacopo Stoppa,   Gabor Szekelyhidi and Richard  Thomas for  helpful discussions. 

\section{Basic definitions}

\subsection{Algebraic theory for group actions}

Let $U,V$ be complex vector spaces of dimensions $q,r$ respectively. We suppose that $G=SL(U)$ acts on $V$ and hence on $\bP=\bP(V)$. We consider the orbit $O$ in $\bP$ of a point $x$. For simplicity we assume that the stabiliser of $x$ in $G$ is finite. We want to study the closure in $\bP$ of $O$. This is an algebraic variety, so for any point $y$ in the closure which is not actually in $O$ we can find a holomorphic map $\Gamma$ from the disc $\Delta \subset \bC$ to $\bP$ with $\Gamma(0)=y$ and $\Gamma(t)\in O$ for $t\neq 0$. We will call such a map an {\it arc through} $y$. Really we should work with germs of such maps, allowing us to restrict to a smaller disc, but we will generally ignore this in our notation. We say that two such arcs $\Gamma_{1}, \Gamma_{2}$ are equivalent if there is a holomorphic map $h:\Delta \rightarrow G$ such that $\Gamma_{1}(t)=h(t) \Gamma_{2}(t)$. Any arc $\Gamma$ can be written in the form 
 $\Gamma(t)= g(t)(x)$ where $g$ is a meromorphic map from $\Delta$ to ${\rm End}(U)$, which restricts to a holomorphic map from the punctured disc $\Delta^{*}$ to $G\subset {\rm End} (U)$. Thus is we choose a vector
$\hat{x}\in V$ lying over $x$ we have a meromorphic map $\hat{\Gamma}$ from $\Delta$ to $V$ defined by $\gamma(t)=g(t) (\hat{x})$.  We define the integer $\nu(\Gamma)$ to be the order of the pole of $\hat{\Gamma}$ at $t=0$. It is obvious that equivalent arcs give the same value of $\nu$.  More explicitly, fixing a basis of $V$ we represent $\hat{\Gamma}$ by its components $(\gamma_{1}, \dots, \gamma_{r})$ say, where $\gamma_{i}$ are meromorphic functions, then  $\nu(\Gamma)$ is the maximum over $i$ of the order of pole of $\gamma_{i}$. Thus $\nu(\Gamma)\leq 0$ if and only if $\hat{\Gamma}$ is holomorphic across
$t=0$ and $\nu(\Gamma)<0$ if and only if $\hat{\Gamma}$ vanishes at $t=0$. Recall that $x$ is called {\it stable} if the orbit of $\hat{x}$ is closed and {\it semistable} if the orbit of $\hat{x}$ does not contain $0$ in its closure.  It is clear that  $x$ is stable if and only if $\nu(\Gamma)>0$ for all $\Gamma$ and semistable if and only if  $\nu(\Gamma)\geq 0$ for all $\Gamma$. 

The most familiar   case of the above set-up is when $g(t)$ has the form
$g(t)= \Lambda(t) R$ where $R$ is a fixed element of $G$ and $\Lambda$ is a {\it one parameter subgroup} in $G$. We call an arc of this type an
{\it equivariant arc}. Replacing $x$ by $Rx$ we can usually reduce to the case when $R=1$. The {\it Hilbert-Mumford criterion} asserts that to test stability or semistability it suffices to restrict to  equivariant arcs. But it is not true that any point $y\in \overline{O}\setminus O$ is contained in an equivariant arc--there are cases where there are boundary points $y$ which are not \lq\lq accessible by 1-parameter subgroups'', see (5.1). This is related to the stabiliser ${\rm Stab}(y)$ of $y$ in $G$. We have
\begin{prop}
\begin{enumerate} \item If ${\rm Stab} ( y) $ is isomorphic to $\bC^{*}$ then any arc through $y$ is equivalent to an equivariant arc.
 \item If ${\rm Stab} (y)$ is reductive then there is an equivariant arc through $y$.
\end{enumerate}
\end{prop}
See (5.2) for the proof.

\

Recall that a {\it weighted flag} in $U$ is a chain of subspaces $U_{1}\subset \dots U_{s}=U$ and associated integers $\lambda_{1}<\lambda_{2}\dots  <\lambda_{s}$. We say that an endomorphism $A$ of $U$ is compatible with the weighted  flag if it preserves the subspaces $U_{i}$ and acts as $\lambda_{i}$ on $U_{i}/U_{i-1}$. It follows that $A$ is diagonalisable, with eigenvalues  $\lambda_{i}$. Conversely, starting with any diagonalisable $A$, with integer eigenvalues, we can define a weighted flag. 

\begin{prop}
Let $g$ be a meromorphic function on the disc with values in ${\rm End}(U)$ holomorphic away from $0$ and with $g(t)$ invertible for $t\neq 0$. Then there is a unique weighted flag in $U$ with the following property. If $A$ is any endomorphism compatible with the weighted flag then we can write
$$  g(t)= L(t) \Lambda(t) R(t)$$
where $L,R$ are holomorphic across $0$ with $L(0), R(0)$ invertible and where $\Lambda(t)= t^{A}$ is the $1$-parameter subgroup generated by $A$.
\end{prop}

This is a standard result. In different language, we can view $g$ as a meromorphic trivialisation over $\Delta^{*}$ of a bundle (in fact the trivial bundle) over $\Delta$. The statement is that in such a situation we get a {\it parabolic structure} on the bundle at the point $0$. To prove the Proposition we can reduce to the case when $g$ is holomorphic (multiplying by a power of $t$). So we can think of it as a map of sheaves
$$   g: \cO^{m} \rightarrow \cO^{m}, $$
and the cokernel of $g$ is a torsion sheaf supported at $0$. Then the statement follows from the classification of torsion modules over $\bC[t]$.

Now we can define  numerical invariants of equivalence classes of arcs $\Gamma$  from the \lq\lq eigenvalues'' $\lambda_{i}$ and their \lq\lq multiplicities'' ${\rm dim} U_{i}/U_{i-1}$ (Of course they are not intrinsically eigenvalues, since the choice of $A$ is not unique, but they are well-defined integers associated to $\Gamma$.) For example we can define the \lq\lq trace''
$$  \sum_{i} \lambda_{i} \ {\rm dim} U_{i}/U_{i-1}.
$$  
This vanishes in the case when $g$ maps into $SL(U)$. We define the \lq\lq norm'' by
$$  \vert \Gamma \vert = \max_{i} \vert \lambda_{i} \vert , $$
which is clearly strictly positive. Now define
\begin{equation} \psi(\Gamma)= \frac{\nu(\Gamma)}{\vert \Gamma\vert}. \end{equation}
This  has the property that it is unchanged if we replace the parameter $t$ by a positive power of $t$.  We can think, informally,  of the numbers $\psi(\Gamma)$, for different arcs $\Gamma$, as a kind of measure of \lq\lq how stable'' the point $x$ is. Next we define
\begin{equation}  \Psi(y) = \sup_{\Gamma} \psi(\Gamma), \end{equation}
where $\Gamma$ runs over the arcs through $y$. In sum, we have attached a numerical invariant $\Psi(y)$ (which {\it a priori} could be $+\infty$) to each point  $y$ in $\overline{O}\setminus O$. 
 
 We have
 \begin{prop}
 The point $x$ is stable if and only if $\Psi(y)>0$ for all points $y\in \overline{O}\setminus O$.
\end{prop}

This is less obvious than it may appear at first sight. The problem is that
we might have a point $y$ such that $\psi(\Gamma)\leq 0$ for some arc through $y$ but not for all. So that while $y$ \lq\lq destabilises'' $x$ we do not have $\Psi(y)\leq 0$. We give the proof in (5.3).

So far our discussion has been entirely algebraic. Now we introduce \lq\lq metric'' geometry. We suppose that we have a hermitian metric on $U$, and hence a maximal compact subgroup $SU(q)\subset G$, and that we have a norm on the restriction of the tautological line bundle over $\bP$  to $\overline{O}\subset \bP$. Since $\overline{O}$ is usually singular it is not quite obvious what we mean by a norm and we leave this point for the moment: for example we could think for the time being of the case when this norm extends smoothly to $\bP$.
We suppose that this bundle norm is invariant under the action of $SU(m)\subset G$. Then we have a moment map
$$  M : \overline{O}\rightarrow \frak{s}\frak{u}(q)^{*}.$$
Now suppose we have $y\in \overline{O}\setminus O$ and $\Gamma$ as above. There is a unique choice of endomorphism $A$, compatible with the weighted flag, which is self-adjoint with respect to the metric on $U$. Then $i A$ lies in $\frak{s}\frak{u}(q)$ and  we get a real number $ \langle M(y), iA\rangle$.

\begin{lem}
In this situation $\nu(\Gamma) \leq \langle M(y), iA\rangle$.
\end{lem}

In the case of equivariant arcs equality holds, which is a more standard and easy fact. To prove the Lemma, write $g(t) = L(t) t^{A} R(t)$ as above. Changing the choice of base point $x$ we can reduce to the case when $L(t)$ is identically $1$ and $R(0)=1$. Now write $V=\bigoplus V_{\mu}$ where $A$ acts with weight $-\mu$ on $V_{\mu}$. So $t^{A}$ acts as multiplication by $t^{-\mu}$ on $V_{\mu}$.
 Let $\omu$ be the largest value of $\mu$ such that $\hat{x}$ has a non-zero component  in $V_{\mu}$. 
Now $R(t) \hat(x)$ has components  $x_{\mu}(t)$ say in $V_{\mu}$ which are (vector valued) holomorphic functions of $t$ and $x_{\omu}(0)$ is non-zero by construction. Thus the $V_{\omu}$ component of $t^{A} R(t) \hat{x}$  grows at least as fast as $t^{-\omu}$ as $t\rightarrow 0$. For $\mu<\omu$ the $V_{\mu}$ component of $t^{A} R(t) \hat{x}$ grows strictly slower than $t^{\omu}$ and for $\mu>\mu_{0}$ it grows strictly slower than $t^{-\mu}$. It follows that $\nu\geq \mu_{0}$ and the components $y_{\mu}$ of $\hat{y}$ vanish if $\mu<\omu$. Let $\mu_{1}$ be the smallest value of $\mu$ such that $y_{\mu}$ is not zero, so we know that $\mu_{1}\geq \omu$. 
Then it follows from the above that we must have $\nu\leq \mu_{1}$ (with strict inequality unless $\mu_{1}=\omu$).

Now let $z\in \bP$ be the limit of $t^{-A}y $ as $t\rightarrow 0$, and let $\hat{z}\in V$ be a representative. Then it is clear that $\hat{z}$ is in $V_{\mu_{1}}$. The definition of the moment map implies that 
$$ \langle M(z), i A\rangle = \mu_{1}. $$
On the other hand the definition also implies that $$ \langle M( t^{-A} y, i A\rangle $$ is an decreasing   function of $t$ so
$$   \langle M(y), iA \rangle \geq \langle M(z), iA \rangle = \mu_{1}. $$
Combining with the inequality $\nu\leq \mu_{1}$ from above we have established the Lemma. 

\

In sum, we now have a \lq\lq differential geometric'' way to obtain a bound on the algebro-geometric invariant $\Psi(y)$. With metric structures as above we have
\begin{equation} \Psi(y) \leq \max_{A}\frac{\langle M(y), iA\rangle}{\Vert A\Vert}, \end{equation}
where $\Vert A\Vert$ denotes the usual operator norm of $A$. Of course, it is equivalent to say that $$\Psi(y)\leq \Vert M(y)\Vert_{1}$$ where $\Vert \ \Vert_{1}$ is the dual \lq\lq trace-norm''. 

\subsection{Application to  Chow varieties}

We will apply this theory to the particular case of Chow varieties. Thus we start with a projective space $\bP(U)$ and for  $n,d$ we consider the set $\Ch$ of $n$-dimensional cycles in $\bP(U)$ of degree $d$. It is a fact that this can be embedded as a projective variety in $\bP(V)$ for a certain representation
$V$ of $G=SL(U)$ but we do not need to know the details of this embedding.  We take a projective manifold $X\subset \bP(U)$ which defines a point in the Chow variety and hence a $G$-orbit.  We also have a Hilbert Scheme
$\Hilb$ which parametrises subschemes with the same Hilbert polynomial as $X$ and a $G$-equivariant regular map from $\Hilb$ to $\Ch$ which is a birational isomorphism. (Note that the scheme structure on $\Hilb$---that is to say, infinitesimal deformations which do not extend to actual deformations--- will not be relevant, so it would be better to take about the underlying variety.)

Now suppose that $Y$ is a an algebraic cycle  in the closure of the orbit of $X$ and choose an arc $\Gamma$ as above. Then $\Gamma$ lifts to a map $\tilde{\Gamma}$ to the Hilbert scheme and $\tilde{\Gamma}(0)$ is a scheme $W$ with underlying cycle $Y$.  The lift $\tilde{\Gamma}$ defines a {\it projective degeneration} of $X$: a flat family $\pi:\cX\rightarrow \Delta$, embedded in $\bP(U)\times \Delta$, with fibre over $t$ isomorphic to $X$ for $t\neq 0$ and to $W$ for $t=0$. Conversely, a projective degeneration of $X$ defines an equivalence class of arcs $\Gamma$. An equivariant arc corresponds to an equivariant degeneration, with a $\bC^{*}$-action, also  called a \lq\lq test configuration''.

Next we go on to the metric theory, so we suppose that $U$ has a Hermitian metric. The basic fact is that there is then a natural induced metric on the restriction of the tautological bundle to $\Ch\subset \bP(V)$, and in particular to the closure $\overline{O}$ of our orbit, so we are in the situation considered above. This theory is explained well by Phong and Sturm  in \cite{kn:PS}. All we really need to know is the corresponding  moment map $M:\Ch\rightarrow \frak{s}\frak{u}^{*}$. Let $A$ be a trace-free self-adjoint endomorphism of $U$ and define 
   a function $H=H_{A}$ on $\bP(U)$ by
$$   H = \frac{1}{\vert x \vert^{2}} \langle x,A x\rangle. $$ 

This is the Hamiltonian for the action of the $1$-parameter group $s\mapsto e^{i A s}$ on $\bP(U)$, with respect to the Fubini-Study symplectic form. Now
let $Z$ be an $n$-dimensional algebraic cycle in $\bP$. Then the formula which defines the moment map is 

$$  \langle M(Z), iA\rangle =  \int_{Z} H d\mu  $$
Here integration over $Z$ is defined in the obvious way, using the volume form $d\mu$ induced by the Fubini-Study metric (normalised so that the volume of $Z$ is equal to to its degree). It is convenient to extend this definition to general Hermitian $A$ by decreeing that the moment map vanishes on multiples of the identity. It will also be convenient to introduce a factor, so we define the \lq\lq Chow number'' of $Z,A$ to be
\begin{equation} \Ch(Z,A)=  \frac{1}{\Vol(Z)} \int_{Z} H d\mu - \frac{\Tr A}{\dim U}. \end{equation}
Thus $\Ch(Z,A)$ is the difference between the average of $H$ over $Z$ and the average eigenvalue of $A$. With this explicit formula in place we can go back to the point we skimmed over before,  involving precisely what we mean by a hermitian structure on a line bundle over a singular space. There is no need to produce a general definition, since one can directly check that the argument above applies in our setting. The crucial points are
\begin{itemize}
\item If $Z$ is preserved by the 1-parameter subgroup generated by $A$ then
$\Ch(Z,A)$ is independent of the metric. On the one hand it is equal to the weight of the induced action on the fibre of the line bundle $\cL^{-1}$ over $Z$, and on the on the other hand it has an interpretation in equivariant cohomology (see \cite{kn:D2}.
\item   $\Ch(Z,A)$ is monotone if $Z$ moves under the 1-parameter subgroup generated by $A$ (see (5.4)). 
\end{itemize}

To sum up, we  attach to each algebraic cycle $Y$ in the closure of the $G$-orbit of $X$ (but not in the orbit) an invariant $\Psi(Y)$. The projective variety $X$ is {\it Chow stable} if and only if $\Psi(Y)>0$ for all such $Y$.  If we choose any Hermitian metric on $U$ we get a bound
\begin{equation}   \Psi(Y) \leq \Vol(Y) \max_{A} \frac{ \Ch(Y,A)}{\Vert A \Vert}. \end{equation}

Notice that from the form of the definition---using arcs---it does not really matter whether we talk about algebraic cycles or schemes here. For a scheme $W$ which is in the closure of the orbit of $X$ in the Hilbert scheme we can define $\Psi(W)$ by taking the supremum over degenerations with central fibre $W$ of the same quantity we used in defining $\Psi(Y)$, so $\Psi(W)\leq \Psi(Y)$  for the cycle $Y$ underlying $W$ and in particular the inequality (5) gives a bound on $\Psi(W)$.

\subsection{ K-stability and   \oK-stability}

We will now think slightly more abstractly of a compact complex $n$-manifold $X$ and positive line bundle $L$ over $X$. Suppose that $L^{m}$ is very ample so its sections define an embedding in $\bP(U)$ with $U=H^{0}(X,L^{m})^{*}$. Suppose that we have a degeneration, as considered above. This is a flat family $\cX\rightarrow \Delta$ with $\cX\subset \bP(U)\times \Delta$. The central fibre is a scheme $W\subset \bP(U)$ and the line bundle $\cO(1)$ on $\bP(U)$ is isomorphic to $L^{m}$ on the non-zero fibres. Take a positive integer $p$ such that $\cO(p)$ is also very ample on all fibres of $\cX$ (this is certainly true if $p$ is sufficiently large). The direct image of $\cO(p)$ is a locally free sheaf over $\Delta$ which can thus be trivialised.
If we fix a trivialisation we get an embedding of $\cX$ in $\bP(U_{p})$ where
$$   U_{p} = H^{0}( X, L^{pm})^{*}. $$

Now the fixed degeneration $\cX$ has a sequence of projective embeddings. We will write $\cX_{p}$ when we want to emphasise the difference, and $W_{p}$ for the central fibre.  For each one we have a numerical  invariant $\nu(\cX_{p})$ say. The notion of {\it K-stability} involves the asymptotics of these as $p\rightarrow \infty$. The theory developed in the literature deals with  the case when $\cX$ is an equivariant degeneration (or test configuration), defined by a $\bC^{*}$ action, so we will now restrict to that situation. Then it is known that $p^{-n} \nu(\cX_{p})$ has a limit as $p\rightarrow \infty$ and we define the {\it Futaki invariant} $F(\cX)$ to be this limit. We say that $(X,L)$ is $K$-stable if $F(\cX)>0$ for all  degenerations $\cX$, and for all $m$. It is important to realise that, even if $L$ is  itself very ample over $X$, degenerations may occur at some large multiplicity $m$ which cannot be   realised by any smaller value, such as $m=1$.

\

Let  $\tilde{X}$ be the blow-up of $X$ at a point $x$. For large integers $\gamma$ the line bundle $\gamma L-  [E]$ is ample on $\tilde{X}$, where $E$ is the exceptional divisor. 

\begin{defn}

The polarised manifold $(X,L)$ is \oK-stable if there is a $\gamma_{0}$ such that the blow-up $\tilde{X}$ at any point of $X$, with the polarisation $\gamma L-[E]$ is $K$-stable, for $\gamma\geq \gamma_{0}$.
\end{defn}

{\bf Note} For the purposes of this paper it is probably the same to say that $(\tilde{X}, \gamma L-[E])$ is {\it asymptotically Chow stable}, thus avoiding the notion of $K$-stability.

\section{Birational modifications}

\subsection{Families over the disc}

We begin by considering a flat family $\pi:\cY\rightarrow \Delta$ over the disc. We write $W$ for the central fibre. Suppose we have an embedding of the family in
$\bP(U^{*}) \times \Delta=\bP\times \Delta $ such that each fibre $\pi^{-1}(t)$, for nonzero $t$, maps to a smooth projective variety $V_{t}$. We suppose that for all $p\geq 1$ and all $t\in \Delta^{*}$ the restriction map 
$$  ev_{t}:H^{0}(\bP; \cO(p)) \rightarrow H^{0}(V_{t}, \cO(p)) $$ is surjective.
We  also suppose that the central fibre $W$  contains a component $B$ which is reduced at its generic point. For each $p\geq 1$ we will define another flat family
$\cY'\rightarrow \Delta$ through a certain birational modification of $\cY$.

 Shrinking the disc if necessary we can suppose that the kernel of $ev_{t}$ has a fixed dimension. This family of kernels has a limit as $t$ tends to $0$. Choose a fixed subspace $J\subset H^{0}(\bP,\cO(p))=s^{p}(U)$ which is a complement to this limit. Then the   $ev_{t}$ to $J$ yields an isomorphism for all $t\in \Delta^{*}$ and when $t=0$ the restriction of the map $ ev_{0}:H^{0}(\bP; \cO(p)) \rightarrow H^{0}(W, \cO(p))$ to $J$ has the same image as $ev_{0}$. Then we get an embedding of $\cX$ in $\bP(J^{*})\times \Delta $ which is just the composite of the original embedding, the Veronese map and a linear projection.

\

Now let $s$ be a nonzero element of $J$, so $e_{0}(s)$ is  a section of 
$\cO(p)$ over  the central fibre. We consider an extension $\sigma$ of  $s$ over $\cY$. Thus $\sigma$ will have the form $\sigma= s+ \sum_{i\geq 1} t^{i}\tau_{i}$
for $\tau_{i}$ in $J$. Such an extension has an {\it order of vanishing} $\nu(\sigma)$ on the component $B\subset \cW$. 
\begin{lem}
For each $s\in J\setminus \{0\}$ there is a $\onu(s)$ such that $\nu(\sigma)\leq \onu(s)$ for all extensions $\sigma$ of $s$.
\end{lem}

Given $\mu>1$ the set of $\tau_{1}$ such that $s+\tau_{1} t$ admits a higher order extension  vanishing to order $\mu$ along $B$ is an affine subspace
$K_{\mu}\subset J$. Clearly the $K_{\mu}$ decrease as $\mu$ increases and hence they are eventually constant. So if there is no such upper bound $\onu(s)$ we can find a $\tau_{1}$ such that $s+\tau_{1} t$ admits  higher order extensions vanishing to arbitrarily high order. Repeating the argument we can find a sequence
$\tau_{i}$ such that for each $\mu$ the finite sum $s+\sum_{i=1}^{\mu} \tau_{i} t^{i}$ vanishes to order $\mu$ along $B$. Thus we get a formal power series $s+\sum_{i=1}^{\infty} \tau_{i} t^{i}$ which vanishes to infinite order. Standard general arguments show that this formal power series can be arranged to be convergent and this contradicts the fact that $i_{t}$ is an isomorphism for $t\neq 0$.

Of course we specify $\onu(s)$ by defining it to be the least possible upper bound  and we set $\onu(0)=+\infty$.  Define 
$$J_{\mu}= \{ s \in J: \onu(s)\geq \mu\}. $$
Then the $J_{\mu}$ are linear subspaces of $J$, defining a flag. Choose a corresponding direct sum decomposition
$$  J= \bigoplus I_{\mu}, $$
where $$  J_{\mu'}= \bigoplus_{\mu\geq \mu'}  I_{\mu}. $$
(Here of course the sums run over a finite subset of integers $\mu$.)
Let $M:J^{*}\rightarrow J^{*}$ be the endomorphism which acts as multiplication by $\mu$ on $I_{\mu}^{*}$. Then for non-zero $t$, we have an automorphism $t^{M} $ of $J^{*}$ and hence of $\bP(J^{*})$. Now we set
$$  V'_{t}= t^{M}(V_{t}). $$
This gives a family $\cV\subset \bP(J^{*})\otimes \Delta^{*}$ over the punctured disc and by general theory there is a unique way to extend this to a scheme $\cY'\subset \bP(J^{*})\times \Delta$, flat over $\Delta$. We will see presently that $\cY'$  is independent of the various choices made in the construction. We write $W'$ for the central fibre of $\cY'$.

\

Notice that we can perform this construction with $p=1$ and then $W'=W$ if and only if $B$ does not lie in any hyperplane in $\bP$.

\

{\bf A simple example}

Let $V_{t}\subset \bC\bP^{2}$ be a family of cubics degenerating to the union of a conic $B$ and a line $R$ with two intersection points $X,Y$.  Take $p=2$. Then $J$ is the $6$ dimensional space of quadratic polynomials and our decomposition is $J= I_{1} \oplus I_{0}$ where $I_{1}$ is the one dimensional subspace spanned by the polynomial $P$ defining the conic $B$. In the Veronese embedding of the original family $\cX\subset \bP^{5}\times \Delta$ the component $B$ lies in $\bP^{4}= \bP(\bC^{5})$ and the component $R$ lies in a $\bP^{2}\subset \bP^{5}$. The modified family $\cX'\subset \bP^{5}$ is different. The central fibre has just one component $B'$, which is a rational curve with a double point at $[P]$. The effect of the modification is to collapse the component $R$ and identify the two intersection points
$X,Y$. 

\

Fix a basis $s_{\alpha}$ of $J$ compatible with the direct sum decomposition so $s_{\alpha}$ lies in $I_{\mu(\alpha)}$. Let $\sigma_{\alpha}$ be extensions vanishing to maximal order. By construction $\osigma_{\alpha}= t^{-\mu(\alpha)}\sigma_{\alpha}$
extends holomorphically over the  generic point in $B$. Let $f_{\alpha}$ be the restriction to $B$. This is meromorphic section of $\cO(p)$ over $B$. More precisely if we lift to the normalisation of $B$ the only poles of $f_{\alpha}$ will be on the intersection $D\subset B $ of the component $B$ in $W$ and the other components.   
  
\begin{lem}
The $f_{\alpha}$ are linearly independent.
\end{lem}

Suppose that there is a linear relation $\sum c_{\alpha} f_{\alpha}=0$, for $c_{\alpha}\in \bC$. Let $\lambdamax$ be the largest value of $\lambda(\alpha)$ for terms with $c_{\alpha} \neq 0$. Choose ordering so that the relation is
$$  c_{1} f_{1} + \dots c_{q} f_{q}+\sum_{\beta>q} c_{\beta} f_{\beta}=0$$
where $\lambda(\alpha)=\lambdamax$ for $\alpha\leq q$ and $\lambda(\alpha)>\lambdamax$ for $\alpha>q$. Then 
$$  c_{1} \sigma_{1}+ \dots c_{q}\sigma_{q}+ \sum_{\beta>q} t^{\lambdamax-\lambda(\beta)} c_{\beta} \sigma_{\beta} $$ is an extension of $c_{1}s_{1}+ \dots c_{q} s_{q}$ which
vanishes to order at least $\lambdamax+1$ along $B$. But by construction the $s_{1}, s_{q}$ are linearly independent elements of the space $I_{\lambdamax}$ and all elements of this space vanish to order exactly $\lambdamax$ so we have our contradiction.

 Let $R^{+}_{p}$ denote the vector space of meromorphic sections of $\cO(p)$ over $\cY$, holomorphic away from the central fibre. 
 Thus elements of $R^{+}_{p}$ can be written as semi-finite Laurent series $\sum \rho_{i} t^{i}$ where the sum has only finitely many negative terms and the coefficients $\rho_{i}$ are in $s^{m}(U)$. Let $R_{p}\subset R^{+}_{p}$ be the subspace of sections which extend holomorphically over the generic point of $B\subset \subset \cW$.  We know that the extensions $\osigma_{\alpha}$  lie in $R_{p}$. Further we have
\begin{lem}
The elements of $R_{p}$ are exactly the sums of the form $\sum_{\alpha} a_{\alpha}(t) \osigma_{\alpha}$, for holomorphic functions $a_{\alpha}$ on $\Delta$. 
\end{lem}
It is obvious that such a sum does lie in $R_{p}$; we have to establish the converse. Any element $\sigma$ of $R^{+}_{p}$ can be written in the form
$$ \sigma= \sum a_{\alpha}(t) \osigma_{\alpha}$$ where $a_{\alpha}(t)$ is meromorphic on $\Delta $, holomorphic away from $0$. Suppose that at least one of the $a_{\alpha}(t)$ has a pole.  Let $k$ be the maximal order of pole of the $a_{\alpha}$ and suppose that $a_{1}, \dots a_{r}$ have  poles of order $k$ while $a_{r+1}, \dots$ have  poles of lower order. Then
, along $B$, the meromorphic section $\sigma$ is $t^{-k}(b_{1} f_{1}+\dots+ b_{r} f_{r})$ where $b_{i}$ is the coefficient of $t^{-k}$ in $a_{i}$. By the linear independence of the $f_{i}$ this section $\sigma$ cannot be extended over the generic point of $B$, so does not lie in $R_{m}$.

\

Now suppose we made different choices in our construction (of the subspace
$J$, of the extensions with maximal vanishing order, of the summands $I_{\mu}$).
We get another collection of sections, $\osigma^{*}_{\beta}$ say, and by the preceding Lemma we can write
$$  \osigma^{*}_{\beta}= \sum g_{ \beta\alpha}(t) \osigma_{\alpha}, $$
for a holomorphic matrix-value function $(g_{\alpha \beta}(t))$. Symmetrically the inverse is also holomorphic,  so defines an automorphism of $\bP(J^{*})\times \Delta$, covering the identity on $\Delta$. It follows from the definition that this automorphism takes the family $\cY'$, constructed using one set of choices, to that constructed using the other set of choices. So we have

\begin{cor}
The family $\cY'\subset \bP(J^{*})\times \Delta $ is uniquely defined, up to automorphisms of $\bP(J^{*})\times \Delta$ covering the identity on $\Delta$.
\end{cor}

We record three simple properties of this construction.

\begin{prop}
\begin{enumerate} 
\item The birational map from $\cX$ to $\cX'$ maps $B$ to a component $B'\subset W'$ which it is reduced at its generic point.
\item We have
$$  \frac{\deg(B')}{\deg(W')}\geq \frac{\deg(B)}{\deg(W)}. $$
\item Suppose we start with $\cX'_{p}$ and perform the same construction with an integer $q$, so we obtain a family $(\cX'_{p})'_{q}$ say. Then
$(W'_{p})'_{q}$  isomorphic to $W'_{pq}$.
\end{enumerate}
\end{prop}

The proof is sketched in (5.6). (The author is grateful to Julius Ross for pointing out item (3) above.)
\

Now as $p$ varies we get a collection of flat families, say $\cY'_{p}\rightarrow \Delta$ with central fibres
$W'_{p}$. Although it is not something we need,  one can also discuss these in the language of graded rings. Let  $R_{\cW}=\bigoplus_{p} R_{p}$ (notice that the definition makes sense for any $p\geq 0$). This is a graded ring with the functions on $\Delta$ as the sub-ring $R_{0}$. We can also evaluate elements of $R_{p}$ on $B$: this gives a vector space, $Q_{p}$, say,  of meromorphic functions on $B$ and $Q=\bigoplus Q_{p}$ is a graded ring. 
Now we can ask whether either
\

\begin{itemize}

\item The ring $R$ is finitely generated over $R_{0}$;
\item The ring $Q$ is finitely generated.
\end{itemize}

It seems reasonable to expect that there properties should be related to the question of the {\it stabilisation} of $\cX_{p}$ or $W_{p}$ as $p$ tends to infinity (or perhaps tends to infinity through multiples of some given integer). Since we do not really need to know this, and the author lacks the relevant background knowledge, we do not go into the matter further here (but see (5.5))

\subsection{b-stability}

Let $(X,L)$ be a polarised manifold and fix $m$ such that $L^{m}$ is very ample. We can also choose $m$ large enough so that the sections of $L^{m}$ generate those of $L^{pm}$ for all $p$.  Let $\cX\subset \bP(U)\times \Delta$ (where $U=H^{0}(X,L^{m})^{*}$) be a projective degeneration with central fibre $W$. We say that $\cX$ is \lq\lq admissible''if it contains a component $B$, as considered above, which is \lq\lq large'' in that the degree (i.e. volume)  of $B$ is greater than half the degree of $X$. Then we call $W$ an {\it admissible limit at multiplicity $m$}.  For large enough $p$ we can apply the construction of the previous subsection to get a new degeneration $\cX'_{p}$ with central fibre $W'_{p}$. We call $W'_{p}$ a {\it  descendant of $W$ at the power $p$}.
We define a \lq\lq web of descendants'' at multiplicity $m$ to be a sequence of schemes $(W'_{1}, W'_{2}, W'_{3}, \dots)$ such that
\begin{enumerate}\item For each $p$, the scheme $W'_{p}$ is an admissible limit at multiplicity
$pm$.
\item  For all $p,q$ the scheme $W'_{pq}$ is a descendant at the power $q$ of $W'_{p}$.
\end{enumerate}

({\bf Remark} Note that this definition includes the statement that $W=W'_{1}$ is a descendant of itself, which just means that $B\subset W$ does not lie in any hyperplane.)

Thus by Proposition 4, any family $\cX$, and in particular any test configuration, defines a web of descendants and if ${\rm Aut} W=\bC^{*}$ any web of descendants beginning with $W$ is obtained in this way. We define the \lq\lq birationally modified Futaki invariant'' of a web of descendants to be
$$  F_{b}= \limsup_{p\rightarrow \infty}  p^{1-n} \Psi(W'_{p}). $$
We define a  {\it destabilising object} for $(X,L)$, at multiplicity $m$, to be a web of descendants at multiplicity $m$ with $F_{b}\leq 0$.

\begin{defn}
   The polarised manifold $(X,L)$ is $b$-stable at multiplicity $m$ if there are no destabilising objects for $(X,L)$ at multiplicity $m$.
   We say that $(X,L)$ is b-stable if it is b-stable at multiplicity $m$ for all large enough $m$. 
\end{defn}

\

The upshot of all this is that we can formulate a new precise form of Yau's conjecture: that {\it the existence of a Kahler-Einstein metric should be equivalent to b-stability}. Of course this is in the case of a Fano manifold
$X$ with positive line bundle $L=K_{X}^{-1}$.

\

{\bf Discussion}

1. The definition of b-stability may seem a little complicated, but much of the complication arises from the fact that we take account  of the possibility of different degenerations with the same central fibre, and hence (perhaps) different descendants of 
$W$ at the same power $p$.  In the case of a web of descendants derived from  ${\rm Aut} W=\bC^{*}$  things are much simpler. Then by Proposition 1 we only need to consider a test configuration and its web of descendants. Suppose we are in the case when $B$ is the whole of $W$. Then our birational modification construction becomes vacuous and the descendant $W'_{p}$ is isomorphic to $W$ but embedded by the linear system $\cO(p)$. The $\bC^{*}$-action on the sections of $\cO(p)$ over $W$ has a generator $A_{p}$ and it is a standard fact that $\Vert A_{p}\Vert= p \Vert A_{1}\Vert$. So 
  $$  \Psi(W_{p}')= p^{-1} \Vert A_{1} \Vert^{-1} \nu(W_{p})$$ and
  hence $ p^{1-n} \Psi(W_{p}')$ tends to the  limit $\Vert A_{1} \Vert^{-1} F(\cX)$, as $p\rightarrow \infty$, where $F$ is the usual Futaki invariant. So, up to the fixed positive factor $\Vert A_{1}\Vert$, our definition reproduces the usual Futaki invariant in this case. It is  possible that for any web of descendants $p^{1-n} \Psi(W'_{p})$ has a limit as $p$ tends to infinity, but because we do not know this we take $\limsup$ instead, in the definition.

2. The definition of b-stability is meant to have a   provisional character. Modifying the definition of an \lq\lq admissible'' degeneration allows us to adjust the definition of b-stability by a notch or two, making the notion less stringent or more. For example we could consider imposing a condition that the degree of the complement of $B$ is very small relative to that of $B$, or that the automorphism group of $W$ is reductive. As we explained in the introduction, such adjustments shift the difficulty from one side of the problem
(b-stable $\Longrightarrow$ KE metric) to the other (KE metric$\Longrightarrow$ b-stable).

\subsection{Families over a general base}

So far, we have considered families over a $1$-dimensional base, in fact the disc $\Delta$. Now we want to discuss the general situation of a flat family $\cX$ over a general variety $N$ with a base point which we write as $0\in N$. As usual, we are really working with germs, so we can  shrink the neighbourhood if necessary. We suppose, first, that the fibre $W$ over $0$ contains a component $B$ and, second, that there is a dense open set $N_{0}\subset N$  such that the  fibres over $N_{0}$ are smooth, just as before. It is not possible to directly extend our construction to define a new family $\cX'$ over $N$ but this can be dome after  blowing up $N$ suitably. Thus given a positive integer $p$ we want to construct a variety
$\hat{N}$ with a regular birational isomorphism $q:\hat{N}\rightarrow N$and a family $\cX'_{p}\rightarrow \hat{N}$ which is characterised by the following property. Any map $\gamma:\Delta\rightarrow N$ with $\gamma(0)=0$ lifts to $\hat{\gamma}:\Delta\rightarrow \hat{N}$ and the pull-back by $\hat{\gamma}$ of $\cX'_{p}$ is isomorphic to the family we have constructed in (3.1), beginning with the pull back $\gamma^{*}\cX$ of $\cX$ to $\Delta$.

 This construction of $\hat{\cX}\rightarrow \hat{N}$ is most likely rather routine, as a matter of algebraic geometry, but we will develop it in a way which is adapted to our differential-geometric application in the next section. 
\

We begin with a basic general fact.

\begin{lem}
Let $\Omega$ be a connected complex manifold and let $f_{0}, \dots f_{m}$ be holomorphic functions on $\Omega \times N$. For $\uX\in \bC^{m+1}$ write $f^{\uX}=\sum X_{\alpha} f_{\alpha}$. Given an pre-compact open subset $\Omega_{0}\subset \Omega$ we can find a finite set of points $z_{1}, \dots z_{r}\in \Omega_{0}$, a neighbourhood $A$ of $0$ in $N$ and $C>0$ such that for all $z\in \Omega_{0}$, $\tau\in A$ and  $\uy\in \bC^{m+1}$ we have
$$    \vert s^{\uy}(z,\tau)\vert \leq C \sum_{i=1}^{r} \vert s^{\uy}(z_{i}, \tau) \vert. $$
\end{lem}

To begin note that it suffices to prove this for $[\uX]$ in a small neighbourhood in $\bC\bP^{m}$, using the compactness of projective space. So write $\uX=(1,y_{1}, \dots, y_{m})$ and define a function $S$ on $\bC^{m}\times \Omega \times N$ by
$$S(\uy, z, \tau)= f_{0}(z,\tau)+ \sum y_{i} f_{i}(z,\tau). $$
For each fixed $z\in \Omega$ we get a function $S(\ , z,\ )$ on $\bC^{m}\times N$, let $I$ be the ideal generated by all of these functions. Now restrict to the local ring of germs of functions about $(0,0)\in \bC^{m}\times N$. This is Noetherian so the ideal is generated by a finite number of functions
$\psi_{i}(\uy, \tau)= S(\uy , z_{i}, \tau )$. It is clear that we can choose the $z_{i}$ to lie in any given open set, in particular in $\Omega_{0}$. Now let $\zeta$ be any fixed point of $\Omega_{0}$. Let $J$ be the ideal of functions $g(\uy, z,\tau)$ in the local ring at $(0,\zeta, 0)$  such that for each fixed $z$ the function $g( \ , z,\ )$ lies in $I$. Then $J$ is finitely generated and by Nakayama's lemma we can take the generators to be polynomials in $z$. The co-efficients of these polynomials must lie in $I$ and it follows that any $g\in J$ can be written in the form
$$   g(\uy, z,\tau)= \sum \chi_{i}(\uy,z,\tau) \psi_{i}(\uy,\tau). $$
In particular $S$ lies in $J$, by construction so we can write
$$  S(\uy,z,\tau)= \sum \chi_{i}(\uy, z,\tau) S(\uy, z_{i},\tau), $$
for points $z$ in a small neighbourhood of $\zeta$ and for small $\uy,\tau$.
Then for such points the desired inequality holds, with $C=\max \vert \chi_{i}(\uy,z,\tau)\vert$.
Now use the precompact hypothesis to cover $\Omega_{0}$ with a finite number of such small neighbourhoods.

Now we proceed with our construction.
Just as before we can fix a subspace $J\subset s^{p}(U)$ and an embedding of $\cX$ in $P(J^{*})\times N$. 
Choose an open set $\Omega$, biholomorphic to a polydisc say, whose closure lies in the smooth part of $B$. Then the family $\cX\rightarrow N $ can be trivialised around $\Omega$, in the sense that, after possibly restricting to a smaller neighbourhood of $0\in N$, there is a holomorphic embedding  $\iota:\Omega\times N\rightarrow\cX$ compatible with the projection $\pi:\cX\rightarrow N$. We also fix a trivialisation of $\iota^{*}(\cO(p))$. Then any element $s$ of $J$ defines a holomorphic function $s(z,\tau)$  on $\Omega \times N$. We choose a large finite set $F$ in $\Omega$ as in the Lemma above, adapted to these functions $s(z,\tau)$.  For each $\tau\in N$ we have a linear map  $e_{\tau}: J\rightarrow \bC^{r}$  defined by $e_{\tau}(s)=(s(z_{1},\tau), \dots s(z_{r},\tau))$. When $\tau$ is in the open dense subset $N_{0}$ this map $e_{\tau}$ is  injective, because an element of the kernel has to vanish on all of $\Omega$ by our choice of $F$ and hence on the whole fibre of $\cX$ (since this is irreducible). Thus, taking the images of the $e_{\tau}$, we get a map from $N_{0}$ to the Grassmannian of $m+1$ dimensional subspaces of $\bC^{r}$. This map need not extend to $N$ but taking the closure of the graph, we can find a blow-up $\hat{N}$ containing a copy of $N_{0}$, to which the map does extend. Let $E$ be the pull-back to $\hat{N}$ of the tautological bundle over the Grassmannian. From this point of view the maps $e_{\tau}$ give a trivialisation of $E$ over $N_{0}$. Thus for $\tau\in N_{0}$ we can map the fibre $\pi^{-1}(\tau)\subset \cX$ into $\bP(E^{*})_{\tau}$. So we get a subvariety ${\cX}_{0}'$ in the projective bundle $\bP(E^{*})$ over $N_{0}$. After perhaps blowing up $\hat{N}$ further, the closure of this defines a flat family ${\cX}'$ over $\hat{N}$. 

Suppose that we make this construction in the case when $N$ is the disc $\Delta$, so no blowing up is required and we get a map from $\Delta$ to the Grassmannian.
Clearly this maps $0$ to the subspace of $\bC^{r}$ with basis  $f_{\alpha}(z_{i})$, in the notation of (3.1),  and one sees that this construction agrees with the previous one. Likewise for the case when we pull back ${\cX}'$ by a map from $\Delta$ to $\hat{N}$. 
\

We can make this construction in the case when $N$ is the closure of the orbit of $X$ in the Hilbert scheme, $N_{0}$ is the orbit and $0$ corresponds to some limiting scheme $W$. Then {\it the descendants of $W$ at the power $p$ are exactly the schemes  parametrised by $q^{-1}(0)\subset \hat{N}$}.
For, in one direction, an arc $\Gamma$ through $0$ has a unique lift to $\hat{\Gamma}$ and we take $\hat{\Gamma}(0)$ which is a point in $q^{-1}(0)$. In the other direction, given a point $\hat{\tau}$ in $q^{-1}(0)$ we can find an arc $\tilde{\Gamma}$ through $\hat{\tau}$ and we obtain the corresponding scheme as a descendant by starting with the arc $q\circ \tilde{\Gamma}$ in $N$.

\

{\bf Remark} Suppose that we adjoin an extra point $z_{r+1}$ to $F$. Then it is clear that the family $\cX'$ we construct will be the same, up to isomorphism. It follows easily that the construction does not depend on the choice of the set $F$, the choice of trivialisations etc.

\section{Applications of the definition}

In this section we will   outline  the relevance of the definitions above to the existence problem for Kahler-Einstein metrics. So we consider a Fano manifold $X$ and take $L=K_{X}^{-1}$. For simplicity we assume that the automorphism group of $X$ is finite. 

\subsection{KE metric $\Longrightarrow$ b-stability}

We have
\begin{prop}
If $X$ has a Kahler-Einstein metric then $(X,K_{X}^{-1})$ is \oK-stable.
\end{prop}
This is a consequence of the results of the author \cite{kn:D2}, Arezzo-Pacard \cite{kn:AP} and Stoppa \cite{kn:S}. For, according to Arezzo and Pacard, the blow up of $X$ admits a constant scalar curvature metric when the parameter $\gamma$ is sufficiently large. Then \cite{kn:D2} shows that the blow-up is at least K-semistable and the refinement of Stoppa shows that it is actually K-stable. (The proof of this refinement will involve blowing up a second time.) In short:  KE $\Longrightarrow$ \oK-stable. It seems reasonable to hope that the argument of Stoppa in \cite{kn:S} can be extended to show that 
$$  \overline{K}-{\rm stability} \Longrightarrow {\rm b-stability} , $$
and, assuming this can be done, we get KE$\Longrightarrow$ b-stability.

In fact one might hope ultimately  to prove the chain
$$ {\rm  KE\ metric}\Longrightarrow \overline{K}-{\rm stable}\Longrightarrow {\rm b-stable}\Longrightarrow {\rm KE metric}, $$

and if this could be done  it would be just as good to take a formulation of the main conjecture involving \oK-stability in place of b-stability. This means that the result in one direction (KE $\Longrightarrow$ \oK- stability) is already in place, but exactly the same work is involved in the extra difficulty of proving the converse (\oK-stability $\Longrightarrow$ KE). So it is really a matter of taste which formulation one prefers. The definition of \oK-stability is quicker to state but puts into prominence the blow-up, which we prefer to see as a device used in the proofs, rather than something fundamental to the problem.

 \subsection{b-stability$\Longrightarrow$ KE metric} 
 
 Here we will  discuss a model problem which does not bear immediately on the general existence question but  which  illustrates  ideas which  apply in other, more complicated, situations. 
 We suppose that we have a sequence of Kahler-metrics $\omega_{i}$  on $X$, in the class $c_{1}(X)$ and that $\Ric(\omega_{i})-\omega_{i}$ tends to zero in $C^{\infty}$, in the sense that for all $l\geq 0$ 
  $$  \max_{X} \vert \nabla^{l} (\Ric(\omega_{i}-\omega_{i}) \vert $$
  tends to zero as $i$ tends to infinity. (Here $\nabla^{l}$ denotes the iterated covariant derivative.) We will also assume that $X$ has complex dimension 3, although this is not fundamental. Then we have
\begin{thm}
If $X$ is b-stable then it admits a Kahler-Einstein metric
\end{thm}

This statement is not intended to be optimal (for example, the arguments could probably be made to work if we restrict to $l=0$), but it has some  content. If $X$ is a suitable small deformation of the Mukai-Umemura manifold $X_{0}$, then  Tian showed in \cite{kn:T1} that $X$  does not admit a Kahler-Einstein metric, but since $X_{0}$ has a Kahler-Einstein metric (as explained in \cite{kn:D3}), there is a sequence $\omega_{i}$ on $X$ satisfying the condition above. (See also the related results of Sun and Wang  \cite{kn:SW} in terms of Ricci flow.)

We do not give a complete proof of the Theorem here, but we will give the part of the argument which brings in the b-stability condition, taking as input four \lq\lq Hypotheses''.

Before beginning it may be helpful to emphasise an elementary but important general point. Suppose we have a sequence of projective spaces $\bP_{i}$, all of the same dimension $N$. Suppose we have varieties $V_{i}\subset \bP_{i}$, all of the same degree. Does it make sense to take the \lq\lq limit'' of the $V_{i}$? Certainly we can choose isomorphisms $\chi_{i}$ from $\bP_{i}$ to the standard model $\bC\bP^{N}$ and then (at least after passing to a subsequence) we can take a limit of the $\chi_{i}(V_{i})$. But if we change the isomorphisms by automorphisms $g_{i}$ of $\bC\bP^{m}$ then the limit of the sequence $g_{i}\chi_{i} V_{i}$ may be completely different. So, as the questions stands, the limit has no intrinsic meaning. Suppose now that the $\bP_{i}$ are \lq\lq metrized projective spaces'' {\it i.e.} projectivisations of hermitian vector spaces. Then we can choose the isomorphisms $\chi_{i}$ to preserve metrics, the automorphisms are reduced to the compact group $PU(N)$ and the limits we get are isomorphic. So the answer to the question is positive if we work with metrised projective spaces. Essentially, for the purposes of taking limits, we can treat metrized projective spaces as being canonically isomorphic.

Now we begin the proof. For each $i$ and  all $k>0$ we get a standard $L^{2}$-norm
$ \Vert\ \Vert_{i,k}$  on 
$H^{0}(X,L^{k})$, using the metric $\omega_{i}$. For large enough $k$, these sections give a \lq\lq Tian  embedding'' $T_{k,i}: X\rightarrow \bP^{N_{k}}$ of $X$ in a metrized projective space of dimension $N_{k}$.  Of course, as above,  we can identify this with the standard space $\bC\bP^{N_{k}}$ with the standard metric.  For fixed $k$ and varying $i$ these embeddings differ by the action of $SL(N_{k}+1, \bC)$ so $T_{k,i}= g_{k,i} \circ T_{k,0}$ say. 

\

{\bf Hypothesis 1} We can fix a large $m$ so that if the sequence $g_{m,i}$ has a bounded subsequence, then $X$ admits a Kahler-Einstein metric.

\

 For any self-adjoint endomorphism $A$, with respect to the $L^{2}$ metric,
$A$ we have  a Chow number $\Ch( T_{k,i} X, A)$. 

\

{\bf Hypothesis 2} There is a function $\epsilon(k)$ with $\epsilon(k)\rightarrow 0$ as $k\rightarrow \infty$ such that,  for all $A$ and $k$ we have
$$ \vert \Ch(T_{k,i} X,A)\vert \leq \Vert A\Vert \epsilon(k). $$
 
 \
 
 (We recall that $\Vert A\Vert$ denotes the operator norm: the modulus of the largest eignevalue.)

\

Now focus attention on the case $k=m$. By compactness of the Hilbert scheme we can suppose (taking a subsequence) that the projective varieties $T_{m,i} X$ have a limit $W$
which is a subscheme of $\bC\bP^{N_{m}}$. If $W$ is projectively equivalent to  $T_{p,0} X$ then the $g_{m,i}$ are bounded and we conclude that $X$ has a Kahler-Einstein metric, by Hypothesis 1. So we suppose the contrary, that $W$ is not equivalent to $X$ and we want to show that $X$ is not b-stable.

\

{\bf Hypothesis 3} The scheme $W$ is an admissible degeneration, that is, it contains a component $B$ of degree greater than one half the degree of $W$  which is reduced at its generic point.

\

Now consider the universal family $\cU$ over a neighbourhood $N$ of $0$ (the point corresponding to $W$) in the closure of the orbit of $X$ in the  Hilbert scheme. So  the sequence $T_{p,i} X$ yields  a sequence $\tau_{i}\in N$ converging to $0$. Fix an open set $\Omega\subset B$ as in (3.3), whose closure lies in the reduced, smooth part of $B$. Then (after perhaps shrinking $N$) we can trivialise the universal family $\cU$ in a neighbourhood of $B$ and define an analytic embedding $\Omega \times N\rightarrow \cU$ compatible with the projection $\cU\rightarrow N$. Thus we have an open set $\Omega_{\tau}$ in the fibre of $\cU$ over $\tau$ and in particular  open sets $\Omega_{\tau_{i}}=\Omega_{i}$ in   $ T_{p,i} X$ which, in an obvious sense, tend to $\Omega$ as $i\rightarrow \infty$.  Now we define a norm $ \Vert \ \Vert_{\Omega, i}$ on sections of $L^{pm}$ by restricting to $\Omega_{i}$ and using the $L^{2}$ norm induced by the standard fibre metric and the Fubini-Study volume form.
\

{\bf Hypothesis 4} The $L^{2}$ norm defined by $\omega_{i}$ and the norm
$\Vert \ \Vert_{\Omega,i}$ are uniformly equivalent (i.e. with constants independent of $i$).

\

{\bf Remark} Note that this implies that the exact choice of $\Omega$ is unimportant: any two choices give equivalent norms. Similarly for the choice of trivialisations etc. 
\

Now we are all set up to state:
\begin{prop}
Assuming Hypotheses (1)-(4) above, if $X$ does not admit a Kahler-Einstein metric then it is not b-stable at multiplicity $m$.
\end{prop}

This will be a consequence of the following, which is the central result of this paper. For each $p>0$ we can suppose that the sequence $T_{pm,i} X$ converges to some scheme $W''_{p}$

\begin{prop}
The sequence $W''_{p}$ is a web of descendants. 
\end{prop}

We give the proof of Proposition 6 assuming Proposition 7. For any self-adjoint endomorphism $A$ we have $\Ch(W''_{p},A)\leq \epsilon(pm) \vert A \vert$ and this implies by (5) that $p^{1-n} \Psi(W''_{p})\leq \epsilon(pm) $. Since $\epsilon(pm)$ tends to $0$ as $p$ tends to infinity we see that $F_{b}\leq 0$, so this web of descendants is a destabilising object.   
\

The proof of Proposition 6 is also easy, given the background we have developed.
It suffices to show that each $W''_{p}$ is a descendant of $W$. (For the other conditions for a web of descendants follow by replacing $m$ by a multiple.) The orbit of $X$ yields an open subset $N_{0}$ of $N$. We have a blow up $q:\hat{N}\rightarrow N$ and a lift $L:N_{0}\rightarrow \hat{N}$.
We know that points of $q^{-1}(0)$ correspond to descendants of $W$. Our sequence $\tau_{i}$ is a sequence in $N_{0}$ converging to $0$ and we thus have a sequence $L \tau_{i}$ in $\hat{N}$. Taking a subsequence we can suppose that $L\tau_{i}$ converges in $\hat{N}$ to some descendant $W'_{p}$. What we need to show is that  this limit is isomorphic to the \lq\lq differential geometric'' limit $W''_{p}$.

Now recall that the family $\hat{\cU}$ over $\hat{N}$ is defined by taking the closure of a family over $q^{-1}N_{0}$. For each point $\tau$ in $N_{0}$ we have a map $e_{\tau}: J\rightarrow \bC^{r}$. The standard Hermitian metric on $\bC^{r}$ restricts to a metric on the image of $e_{\tau}$. On the other hand we can identify $J^{*}$ with the space of sections of $L^{pm}$ over the fibre $V_{\tau}$ in $\cU$. In particular we can do all this for the points
$\tau_{i}$ so that $V_{\tau_{i}}$ is $ T_{m,i} X$. Then for each $i$ we have three norms on the space of sections $H^{0}(X,L^{pm})$.
\begin{enumerate}
\item The standard $L^{2}$ norm defined by $\omega_{i}$.
\item The norm $\Vert \ \Vert_{\Omega,i}$ defined by the $L^{2}$ norm over $\Omega$. 
\item The norm defined by the map $e_{\tau_{i}}: J\rightarrow \bC^{r}$, as above.
\end{enumerate}
Unwinding the constructions, to prove that $W'_{p}$ and $W''_{p}$ are isomorphic it suffices to show that the first and third norms are uniformly equivalent, with constants independent of $i$. Hypothesis 4 states that the first norm is equivalent to the second norm, so it suffices to show that the second and third norms are uniformly equivalent. But this is clear from the way the finite set $F$ was chosen (Lemma 5). In one direction, the third norm dominates the $L^{\infty}$ norm of sections over a slightly smaller open set (and we have pointed out above that the precise choice of $\Omega$ is not important). In the other direction the $L^{2}$-norm over $\Omega$ dominates the $L^{\infty}$ norm over an interior set by standard elliptic estimates. 

\

 As we mentioned above, the four Hypotheses hold for the sequence of \lq\lq approximate KE metrics'' $\omega_{i}$, thus proving the Theorem. The proofs will be given elsewhere. Hypotheses 1,3 and 4 are essentially known results, and very similar statements can be found in the recent paper \cite{kn:T2} of Tian. Thus the main new input is Hypothesis 2. The proof of this (which leads us to make the restriction on the dimension of $X$ for the time being) depends in turn on joint work with X-X Chen.

\subsection{Miscellaneous remarks }
\begin{enumerate}
\item {\bf Uniform stability}

A general difficulty which arises in applying the usual definitions of stability is that these definitions state that the Futaki invariant is positive but do not supply any definite lower bound. In the special case of toric manifolds (and in the context of extremal metrics)  establishing such a lower bound was one of the main points of \cite{kn:D1} (Proposition 5.2.2 in \cite{kn:D1}). In general, this issue was considered by Szekelyhidi \cite{kn:Sz} who introduced a notion of \lq\lq uniform stability''. His definition has the shape
$$    F(\cX) \geq c \Vert A \Vert, $$
for all test configurations $\cX$ and 
for an appropriate norm on the generator of the action (and, in a general context, it could be that different norms lead to different notions). Taking such a definition of stability makes it easier to prove the direction \lq\lq stable'' $\Longleftarrow$ KE metric, but it seems hard to establish the converse.
Our definition of b-stability adopts this idea to some extent: roughly it asserts uniform stability over a restricted class of test configurations. The fundamental point is that the result of Arezzo and Pacard gives a little extra control of the Futaki invariants (or Chow weights).

 \
 
 \item {\bf Testing stability}
 
 \
 
 As we mentioned in the Introduction, the impact of this whole discussion is rather limited unless one has a way of testing \lq\lq stability'' in explicit situations. The difficulty (in the case of \oK-stability, say) is that the definition requires checking test configurations of arbitrarily high multiplicity--that is to say, degenerations of $X$ embedded in arbitrarily large projective spaces. For a fixed, reasonably small, multiplicity it may be possible to analyse all the test configurations, but this gets more and more complicated as the dimension grows. In this regard, it is relevant that the argument outlined above in the direction b-stability $\Longrightarrow$ KE metric works with b-stability at an explicit multiplicity $m$, which can be computed in principle from analytical information. If this argument can be refined to produce a multiplicity $m$ which is reasonably small then one could hope to verify b-stability in some explicit cases.

Again, this issue arises in the toric case. In dimension 2, there is a straightforward test (for K-stability) involving certain \lq\lq simple'' piecewise-linear functions on the polytope of the variety (see \cite{kn:D1}). But in higher dimensions, even in the toric case, the situation is less clear.

\item {\bf Finite generation and Gromov-Haussdorf limits}

Recall that we associated a graded ring $Q$ to an admissible degeneration, and we expect that the finite generation of this ring should be related to the stabilisation of the sequence of descendants. It also seems likely that  $Q$ can    be obtained as the limits of $L^{2}$ sections over the Gromov-Hausdorff limit considered by Ding and Tian \cite{kn:DT}. The finite generation question could then  be seen as the problem of endowing this Gromov-Haussdorf limit with an algebraic structure. If $Q$ is finitely generated then the obvious candidate is ${\rm Proj}(Q)$.

\item {\bf Extremal and constant scalar curvature metrics}

The definition of b-stability is aimed at the Kahler-Einstein case and is not intended to be appropriate, as it stands, in the more general setting of extremal and constant scalar curvature metrics. On the other hand an example of Apostolov {\it et al}  \cite{kn:ACGT} shows that K-stability is probably not the correct criterion for the existence of these metrics. It may be the same general idea---involving not just one test configuration but a sequence---will be relevant in that case too.

\end{enumerate}

\section{Examples and subsidiary results}

\subsection{Points not accessible by one-parameter subgroups}

We consider the action of $SL(3,\bC)$ on $s^{d}(\bC^{3})$, polynomials in $x,y, z$,  so points of the projectivisation corresponds to plane curves. Suppose that a curve $C$ meets the line $z=0$ in $d$ points $p_{1},\dots, p_{d}$. Then the action of the 1-parameter group $x\mapsto t x, y \mapsto t y , z\mapsto t^{-2} z$ deforms $C$ into the union of the $d$ lines $Op_{i}$, where $O$ is $x=y=0$. If $d>5$ then simple dimension counting shows that a typical singular curve consisting of $d$ lines through a point is not projectively equivalent to one which arises in this way from the fixed curve $C$. Now consider a curve $C$ of degree $6$ defined by the equation
$$ z \prod_{i=1}^{5} (x-\lambda_{i} y) + p(x,y)=0. $$
Thus there are 5 branches of the curve passing through $x=y=0$.
The inverse of the $1$-parameter subgroup above deforms this curve into the union of the 5 lines $x=\lambda_{i} y$ through $O$ and the line at infinity. Making a projective transformation, fixing lines through $O$, we can move the line at infinity to a line $x-\mu y=z$ say. Now apply the same $1$-parameter subgroup to deform this to the curve $C'$ which is the union of 6 lines through $O$: the 5 lines $x=\lambda_{i} y$ and the sixth line $x=\mu y$. Thus $C'$ is in the closure of the orbit of $C$ but there is no reason why this configuration of $6$ lines should occur from the intersection of $C$ with a line,  so that in general
$C'$ will not be accessible by a $1$-parameter subgroup. 

\subsection{Proof of Proposition 1}
The result is similar to Luna's slice theorem, but the author has not found this exact statement in the literature. For $\xi\in \frak{g}$ and $z$ in $\bP(V)$ we write $\xi z$ for the corresponding tangent vector to $\bP(V)$ at $z$.
We suppose that the stabiliser of $y$ is the complexification $K^{c}$ of a compact group $K$. (The point of the condition is that  representations of $K^{c}$ decompose as sums of irreducibles.) Then we  find an equivariant slice for the action  a projective subspace $\bP'\subset \bP(V)$ with the properties
\begin{enumerate}
\item $\bP'$ contains $y$ and  is preserved by $K^{c}$;
\item $\bP'$ is transverse  to the $G$-orbit of $y$ at $y$;
\item for $y'\in \bP'$ near $y$ and $\xi \in \frak{g}$ the tangent vector $\xi y'$ lies in the tangent space of $\bP'$ if and only if $\xi$ is in $\frak{k}$.
\end{enumerate}

  To do this we decompose $\frak{g}=\frak{k}^{c}\oplus \frak{p}$ as representations of $K^{c}$. The derivative of the action at $x$ gives us a decomposition of $K^{c}$-representations $V= \bC \hat{x}\oplus \frak{p} \oplus S$ say. Then we can take $\bP'=\bP(\bC\hat{x}\oplus S)$. 
Now the third condition above implies that the identity component of the stabiliser in $G$ of any point in $\bP'$ close to $y$ is a subgroup of $K^{c}$.
Let $V_{0}$ be the intersection of the  of the G-orbit of $O$ of $x$ with $\bP'$.  The second condition implies that the $G$-orbit of any point in $\bP(V)$ close to $y$ meets $\bP'$ in a point which is also close to $y$.
It follows that $y$ lies in the closure $V$ of $V_{0}$.   The third condition implies that if $y'$ is close to $y$ and $y'$ is in $V_{0}$ then the for small ball $B$ around $y'$ the intersection of $B$ with $V_{0}$ is equal to the intersection of $B$ with the $K^{c}$-orbit of $y'$. Since $V$ is an algebraic variety it is clear that, near to $y$, $V_{0}$ is contained in a finite union of $K^{c}$ orbits. So there is a single $K^{c}$ orbit $K^{c} y'$ in $V_{0}$
which contains $y$ in its closure. But the action of $K^{c}$ on $\bP'$ is determined by the linear action on $S$ and it follows from the Hilbert-Mumford criterion, applied to the reductive group $K^{c}$, shows that there is a $1$-parameter subgroup $\Lambda$ of $K^{c}$ such that $y$ lies in the closure of the $\Lambda$-orbit of $y'$. Viewing  $\Lambda$ as a 1-parameter subgroup in $G$ we obtain the first statement of the Proposition. 

For the second statement, we see from the transversality condition (2) that any arc through $y\in \overline{O}$ is equivalent to an arc mapping into $V$. Near to $y$, we know that $V_{0}$ is a finite union of $\bC^{*}$ orbits so the arc must map into a single one of these orbits and it follows that the arc is equivalent to an equivariant arc.

{\bf Remark}
It seems likely that in fact $V_{0}$ is, near to $y$, equal to a single
$K^{c}$-orbit but the author has not managed to prove this. If this were the case it would follow that when the stabiliser is $\bC^{*}$ there is a unique equivalence class of arcs through $y$, up to the obvious fact that we can take a covering, replacing $t$ by a power of $t$. 

\subsection{Proof of Proposition 3}
The argument is related to ideas of Thaddeus \cite{kn:Th}. In one direction, it is clear that if $x$ is stable then $\Psi(y)>0$ for all $y$. In the other direction,
suppose that $x\in \bP$ is a point which is not stable for the $G$-action. We choose some other representation of $G$ and hence another projective space $Q$ on which $G$ acts. We can choose $Q$  to contain a point $q$ which {\it is} stable. We fix metrics so that we have moment maps $\mu_{\bP}, \mu_{Q}$. Now consider the $G$-action on the product $\bP\times Q$.  The notion of stability depends on the choice of a  class in $H^{2}(\bP\times Q)$ which we take to be $\omega_{\bP} +\eta \omega_{Q}$ where $\eta$ is a real parameter, and $\omega_{\bP},\omega_{Q}$ are the standard generators. When $\eta$ is irrational we move outside the algebro-geometric framework, but we can still apply the symplectic theory of Kirwan \cite{K}. The moment map for the product is $\mu_{\bP}+\eta \mu_{Q}$. When $\eta$ is very large the point $(x,q)$ is stable and there is a point in its orbit where $\mu_{\bP}+\eta \mu_{Q}$ vanishes. But when $\eta=0$ there is no point in the orbit where
$\mu_{\bP}$ vanishes, since $x$ is not stable. It follows that there is some $\eta_{0}\geq 0$ and a point $(y,r)$ in the closure of the $G$-orbit in
$\bP\times Q$ such that $\mu_{\bP}(y)+\eta_{0}\mu_{Q}(r)=0$. We claim that
$\Psi(y)\leq 0$. Let  $\Gamma_{\bP}:\Delta\rightarrow \bP$ be an arc through $y$. We can lift this to the $G$-orbit in $\bP\times Q$, so we get a map $(\Gamma_{\bP},\Gamma_{Q}):\Delta\rightarrow \bP\times Q$. The maps $\Gamma_{\bP}, \Gamma_{Q}$ define the same weighted flag and hence the same self-adjoint endomorphism $A$.  We have integers $\nu_{\bP}(y), \nu_{Q}(r)$.  The inequality of Lemma 1 gives $$  \nu_{\bP}(y)\leq \langle \mu_{\bP} (y), iA\rangle\ \ ,\  \nu_{Q}(r)\leq \langle\mu_{Q}(r), iA\rangle, $$
so $ \nu_{\bP}(y)+\eta_{0}\nu_{Q}(r)\leq 0$. But the condition that $q$ is stable implies that
$\nu_{Q}(r)>0$ so we see that $\nu_{\bP}\leq 0$. This implies that $\Psi(y)\leq 0$, as required. 
\subsection{Monotonicity of the Chow invariant}
Let $A$ be a self-adjoint endomorphism of a Hermitian vector space $U$ and let $Z$ be an algebraic cycle in $\bP(U)$. For real $s$ let $Z_{s}= e^{As} Z$ and consider the function $f(s)=\Ch(Z_{s}, A)$. We want to show that $f$ is an increasing function of $s$. This is a well-known fact, first proved by Zhang [Z], and there are a number of proofs in the literature, but we include a short proof here for completeness. Consider for the moment a more general situation of a compact Riemannian manifold $P$, a real valued function $h$ on $P$ and a submanifold
$V$ of $P$. Let $\phi_{s}:P\rightarrow P$ be the gradient flow of $h$ and let $V_{s}=\phi_{s}(V)$. We consider the function
\begin{equation}  g(s)= \int_{V_{s}} h\  d\mu, \end{equation}
where $d\mu$ is the induced Riemannian volume element. Then
$$  g'(s)= \int_{V_{s}} \vert ({\rm grad}_{\perp} h) \vert^{2}d\mu + \int_{V_{s}} \langle M, {\rm grad}_{\perp}h \rangle d\mu, $$
where ${\rm grad}_{\perp} h$ is the component of the gradient vector field normal to $V_{s}$ and $M$ is the mean curvature vector of $V_{s}$. Here the second term arises when we differentiate the volume element in (7). The relevance of this is that when $P$ is complex projective space the flow $e^{As}$ is the gradient flow of the function $H$. Since the mean curvature of a complex subvariety vanishes, the second term drops out and we derive the desired monotonicity.
(It is not hard to see that the presence of  singularities of $Z$ does not affect the argument.)

\subsection{Zariski's example}

We explain how a well-known example of Zariski [Z], fits into the framework of this paper.
  Let $S$ be the blow-up of the projective plane in 12 points $q_{i}$ and $L$ be the line bundle defined by the divisor $D=2( 4H-\sum E_{i})$ on $S$. If the points are in general position then $L$ is very ample on $S$ but if the points lie on a cubic $C\subset \bP^{2}$ then the proper transform $\tilde{C}$ has intersection number $0$ with $D$ so $L$ is not ample. If $q_{i}$ are in general position on $C$ then no multiple $L^{p}$ is trivial on $\tilde{C}$.  All sections of $\cO(pD)$ vanish on $C$ but, for all $p$, there are sections which vanish with multiplicity $1$. It follows that the ring $\bigoplus H^{0}(S;L^{p})$ is not finitely generated.
Now vary $q_{i}$ in curves $q_{i}(t)$ for $t\in \Delta$ and hence construct a family $\cS\rightarrow \Delta$. We suppose that the $q_{i}(t)$ are in general position for non-zero $t$ but $q_{i}(0)$ are in general position on a cubic curve, as above.  We have a line bundle $\cL\rightarrow \cS$ and the direct image is locally free so we can regard the fibres $S_{t}$, for non-zero $t$ as projective varieties in a fixed projective space. Taking the closure we get a family $\cX\rightarrow \Delta$ where the central fibre has two components. One is  $S_{0}$, embedded by sections of $L^{p}\otimes[-\tilde{C}]$, and the other is a $\bP^{1}$ bundle $R$ over $\tilde{C}$. For any $p>0$ we can apply our construction from (3.1) (taking $B=S_{0}$) to form $\cX'_{p}$. This is exactly the same as the family we get by replacing $L$ by $L^{p}$ in the construction of $\cX$ above.  For all $p$ we get the same central fibre $S_{0}\cup_{\tilde{C}} R$  but with a different line bundle over it. The ratio  $ \deg(R)/\deg(S_{0})$ is$O(p^{-1})$ as $p\rightarrow \infty$. The graded ring $Q$ defined in (3.1) is isomorphic to $\bigoplus H^{0}(S_{0}, L^{p})$ and is not finitely generated. 

{\bf Remark} In this example there is a natural family, independent of $p$, in the background--the family $\cS$. The problem is that we are trying to embed this using sections of a line bundle which is not positive on the central fibre, and this gives rise to the infinite series of \lq\lq descendants''.  Whether an example like this can occur for a degeneration of a Fano manifold and, conversely, whether this is the only mechanism by which we can obtain an infinite sequence of \lq\lq descendants'' are questions  beyond the authors knowledge, but which seem to be important.  

\subsection{Proof of Proposition 4}

We can regard $W$ as being embedded in $\bP(J^{*})$. If $B\subset W$ does not lie a hyperplane in $\bP(J^{*})$ then $W'=W$ and we are done, so suppose that there $J^{*}=U_{1}\oplus U_{2}$ and $B$ lies in $\bP(U_{1})$. Let $B_{0}\subset B$ be the complement of the intersection of $B$ with the other components of $W$ and of the support of nilpotents. Then, following through the  construction, we see that $B_{0}$ maps isomorphically to $B'_{0}\subset B'$ by a map of the form $x\mapsto (x,f(x))\in \bP(U_{1}\oplus U_{2})$ where $f$ is a holomorphic section of $\cO(1)\otimes U_{2}$ over $B_{0}$. Thus it is clear that $B'$ is reduced at its generic point. We can choose a linear subspace $\bP(R)\subset\bP(U_{1})$ of complementary dimension such that $B\cap \bP(R)$ lies in $B_{0}$. Then we get the same number of intersection points of the image of $B_{0}$ with $\bP(R\oplus U_{2})$. Since any other intersection points give a positive contribution to the intersection number we see that $\frac{\deg(B')}{\deg(W')}\geq \frac{\deg(B)}{\deg(W)}$. 

To see that $(W'_{p})'_{q}= W'_{pq}$ we can use the approach through the norms on sections of $\cO(pq)$ defined by restricting to a large finite set $F$. Both $(W'_{p})'_{q}$ and $ W'_{pq}$ can be obtained as the limits of the nonzero fibres $V_{t}$ under embeddings in $H^{0}(V_{t}, \cO(pq))^{*}$, so, as in Section 4, we have to show that two norms on this space are uniformly equivalent (with constants independent of $t$). But this is rather clear from the construction using the fact that our local trivialisation of the family $\cX$ around $\Omega\subset B$ maps holomorphically to a local trivialisation
of $\cX'_{p}$ around the image $\Omega'$ by a map of the same form as above.  


\begin{thebibliography}{99}
\bibitem{kn:ACGT} Apostolov, V., Calderbank, D., Gauduchon, P. and Tonnesen-Friedman, C. {\em Hamiltonian 2-forms in Kahler geometry III: extremal metrics and stability} Inventiones  Math. 173 547-60 (2008)
\bibitem{kn:AP} Arezzo, C. and Pacard, F. {\em Blowing up and desingularising constant scalar curvature Kahler metrics, I} Acta Math. 196 179-228 (2006)
\bibitem{kn:DT} Ding, W. and Tian, G. {\em Kahler-Einstein metrics and generalised Futaki invariants} Inventiones Math. 110 315-335 (1992)
\bibitem{kn:D1} Donaldson, S. {\em Scalar curvature and stability of toric varieties} Jour. Differential Geometry 62 289-349 (2002)
\bibitem{kn:D2} Donaldson, S. {\em Lower bounds on the Calabi functional} Jour. Differential Geometry 70 453-472 (2005)
\bibitem{kn:D3} Donaldson, S. {\em Kahler geometry on toric manifolds, and some other manifolds with large symmetry} In: handbook of Geometric Analysis, Vol. I Ji et al eds. International Press 29-74 (2008)
\bibitem{kn:K} Kirwan,F.C {\em The cohomology of quotients in symplectic and algebraic geometry} Princeton UP 1983
\bibitem{kn:PS} Phong, D.H.  and Sturm, J. {\em Scalar curvature, moment maps and the Deligne pairing} Amer. J. Math. 126 693-712 (2004)
\bibitem{kn:RT} Ross, J. and Thomas, R. {\em  A study of the Hilbert-Mumford criterion for the stability of projective varieties} Jour. Algebraic Geometry 10 201-255 (2007)
\bibitem{kn:S} Stoppa, J. {\em K-stability of constant scalar curvature Kahler manifolds} Advances Math. 221 1397-1408 (2009)
\bibitem{kn:SW} Sun, S.  and Wang, Y. {\em On the Kahler-Ricci flow near a Kahler-Einstein metric}  Arxiv 10004.2018
\bibitem{kn:Sz} Szekelyhidi, G. {\em Extremal metrics and K-stability} PhD. Thesis 2006 DG/061102
\bibitem{kn:Th} Thaddeus, M. {\em Geometric Invariant Theory and flips}
Jour. Amer. Math. Soc. 9 691-723 (1996)
\bibitem{kn:T1} Tian, G. {\em Kahler-Einstein metrics of positive scalar curvature} Inventiones Math. 130 1-57 (1997)
\bibitem{kn:T2} Tian, G. {\em Existence of Einstein metrics on Fano manifolds} Preprint 2010.
\bibitem{kn:Y} Yau, S-T. {\em Review of Kahler-Einstein metrics in algebraic geometry} Israel Math. Conference Proc., Bar-Ilan Univ.  9 433-443 (1996)
\bibitem{kn:Z} Zhang, S. {\em Heights and reductions of semi-stable varieties}
Compositio Math. 104 77-105 (1996)
\bibitem{kn:Z} Zariski,O. {\em The theorem of Riemann-Roch for high multiples of an effective divisor on an algebraic surface} Ann. Math. 76 560-615 (1962)

\end{thebibliography}
\end{document}